\documentclass[english,11pt]{article}
\usepackage[latin1]{inputenc}
\usepackage{amsmath,amsthm,amssymb,babel}

\newtheorem{teo}{Theorem}

\newtheorem{lemma}{Lemma}
\newtheorem{rem}{Remark}

\def\proof{{\it Proof.}\ }
\def\endproof{\hfill $\Box$\par\vskip3mm}

\def\eq#1{(\ref{#1})}

\def\neweq#1{\begin{equation}\label{#1}}
\def\endeq{\end{equation}}

\def\phi{\varphi}
\def\RR{{\mathbb R} }

\def\ZZ{{\mathbb Z} }

\def\di{\displaystyle}

\date{}

\title{\sc Nonlinear eigenvalue problems in Sobolev spaces with variable exponent}
\author{Teodora-Liliana Dinu \\ \small Department of Mathematics, ``Fra\c tii Buze\c sti" College, 
Bd. \c{S}tirbei--Vod\u a No. 5, 200352 Craiova, Romania\\ \small E-mail: {\tt tldinu@gmail.com}}

\begin{document}
\baselineskip16pt
\maketitle
\noindent{\small{\sc Abstract}.
We study the boundary value problem $-{\rm div}((|\nabla u|^{p_1(x)
-2}+|\nabla u|^{p_2(x)-2})\nabla u)=f(x,u)$ in $\Omega$,
$u=0$ on $\partial\Omega$, where $\Omega$ is a smooth bounded domain in
$\RR^N$. We focus on the cases when
$f_\pm (x,u)=\pm(-\lambda|u|^{m(x)-2}u+|u|^{q(x)-2}u)$, where $m(x):=\max\{
p_1(x),p_2(x)\}<q(x)<
\frac{N\cdot m(x)}{N-m(x)}$ for any $x\in\overline\Omega$. 
In the first case we show the existence of infinitely many weak 
solutions for any $\lambda>0$. In the second case we prove that if 
$\lambda$ is large enough then there exists a nontrivial weak 
solution. Our approach relies on the variable exponent theory of 
generalized Lebesgue-Sobolev spaces, combined with a 
$\ZZ_2$-symmetric version for even functionals of the Mountain Pass 
Lemma and some adequate variational methods.  
}

\section{Introduction and preliminary results}
Electrorheological fluids (sometimes referred to as ``smart fluids"), 
are particular fluids of high technological interest whose apparent
viscosity changes reversibly in response to an electric field. The
electrorheological fluids have been intensively studied from the 
1940's to the present. The first major 
discovery
on electrorheological fluids is due to Willis M. Winslow \cite{WW}. 
He noticed that such fluids' (for instance lithium
polymetachrylate) viscosity in an electrical field is inversely
proportional to the strength of the field. The field induces 
string-like
formations in the fluid, which are parallel to the field. They can 
raise
the viscosity by as much as five orders of magnitude. This phenomenon 
is
known as the Winslow effect. For a
general account of the underlying physics confer \cite{hal} and for 
some
technical applications
\cite{pfe}.
We just remember that any device which currently depends upon 
hydraulics, hydrodynamics or hydrostatics can benefit from
electrorheological fluids' properties. Consequently, electrorheological
fluids are most promising in aircraft and aerospace applications. 
For more information on properties and the
application of these fluids we refer to
\cite{AM,D,hal,R}.
   
The mathematical modelling of electrorheological fluids determined 
the study of variable Lebesgue and Sobolev spaces $L^{p(x)}$ and
$W^{1,p(x)}$, where $p(x)$ is a real-valued function. Variable
exponent Lebesgue spaces appeared in the literature for the 
first time already in a 1931 article by W.~Orlicz \cite{orl}. In
the years 1950 this study was carried on by Nakano \cite{nak} who
made the first systematic study of spaces with variable exponent. 
Later, the Polish mathematicians investigated the modular function 
spaces (see, e.g., the basic monograph Musielak \cite{M}). Variable exponent Lebesgue 
spaces on the real line have been independently developed by Russian
researchers. In that context we refer to the work of Tsenov   
\cite{tse}, Sharapudinov \cite{sha} and Zhikov \cite{Z1, Z2}. For deep results
in weighted Sobolev spaces with applications to partial differential equations we refer to the 
excellent monographs by  Drabek, Kufner and Nicolosi \cite{dkn}, by Hyers,  Isac and Rassias
\cite{isac}, and by Kufner and Persson \cite{kp}.

Our main purpose is to study the boundary value problem
\begin{equation}\label{Pr1}
\left\{\begin{array}{lll}
-{\rm div}((|\nabla u|^{p_1(x)-2}+|\nabla u|^{p_2(x)-2})\nabla u)=
f(x,u), &\mbox{for}& x\in\Omega\\
u=0, &\mbox{for}& x\in\partial\Omega
\end{array}\right.
\end{equation}    
where $\Omega\subset\RR^N$ ($N\geq 3$) is a bounded domain with
smooth boundary and $1<p_i(x)$, $p_i(x)\in C(\overline\Omega)$
for $i\in\{1,2\}$. We are looking for nontrivial weak
solutions of Problem \eq{Pr1} in the generalized Sobolev space
$W^{1,m(x)}(\Omega)$, where $m(x)=\max\{p_1(x),p_2(x)\}$ for any
$x\in\overline\Omega$. We point out that problems of type \eq{Pr1} were 
intensively studied in the past decades. We refer to \cite{CF,
FZh, FZZ} for some interesting results.  

We recall in what follows some definitions and basic properties
of the generalized Lebesgue--Sobolev spaces $L^{p(x)}(\Omega)$
and $W_0^{1,p(x)}(\Omega)$, where $\Omega$ is
a bounded domain in $\RR^N$.

Set
$$C_+(\overline\Omega)=\{h;\;h\in C(\overline\Omega),\;h(x)>1\;{\rm 
for}\;
{\rm all}\;x\in\overline\Omega\}.$$
For any $h\in C_+(\overline\Omega)$ we define
$$h^+=\sup_{x\in\Omega}h(x)\qquad\mbox{and}\qquad h^-=
\inf_{x\in\Omega}h(x).$$
For any $p(x)\in C_+(\overline\Omega)$, we define the variable exponent
Lebesgue space
$$L^{p(x)}(\Omega)=\{u;\ u\ \mbox{is a
 measurable real-valued function such that }
\int_\Omega|u(x)|^{p(x)}\;dx<\infty\}.$$
We define a norm, the so-called {\it Luxemburg norm}, on this space by 
the
formula
$$|u|_{p(x)}=\inf\left\{\mu>0;\;\int_\Omega\left|
\frac{u(x)}{\mu}\right|^{p(x)}\;dx\leq 1\right\}.$$
Variable exponent Lebesgue spaces resemble classical Lebesgue spaces
in many respects: they are Banach spaces \cite[Theorem 2.5]{KR}, the
H\"older inequality holds
\cite[Theorem 2.1]{KR}, they are reflexive if and only if $1 < p^-\leq
p^+<\infty$
\cite[Corollary 2.7]{KR} and continuous functions are dense if $p^+
<\infty$ \cite[Theorem 2.11]{KR}. The inclusion between
Lebesgue spaces also generalizes naturally \cite[Theorem 2.8]{KR}: if 
$0 <
|\Omega|<\infty$
 and $r_1$, $r_2$
are variable exponents so that $r_1(x) \leq r_2(x)$ almost everywhere 
in $\Omega$ then there exists the continuous embedding
$L^{r_2(x)}(\Omega)\hookrightarrow L^{r_1(x)}(\Omega)$, whose norm 
does not exceed $|\Omega|+1$.

We denote by $L^{p^{'}(x)}(\Omega)$ the conjugate space
of $L^{p(x)}(\Omega)$, where $1/p(x)+1/p^{'}(x)=1$. For any
$u\in L^{p(x)}(\Omega)$ and $v\in L^{p^{'}(x)}(\Omega)$ the H\"older
type inequality
\begin{equation}\label{Hol}
\left|\int_\Omega uv\;dx\right|\leq\left(\frac{1}{p^-}+
\frac{1}{{p^{'}}^-}\right)|u|_{p(x)}|v|_{p^{'}(x)}
\end{equation}
holds true.

An important role in manipulating the generalized Lebesgue-Sobolev 
spaces
is played by the {\it modular} of the $L^{p(x)}(\Omega)$ space, which 
is
the mapping
 $\rho_{p(x)}:L^{p(x)}(\Omega)\rightarrow\RR$ defined by
$$\rho_{p(x)}(u)=\int_\Omega|u|^{p(x)}\;dx.$$
If $(u_n)$, $u\in L^{p(x)}(\Omega)$ and $p^+<\infty$ then the 
following relations hold true
\begin{equation}\label{L4}
|u|_{p(x)}>1\;\;\;\Rightarrow\;\;\;|u|_{p(x)}^{p^-}\leq\rho_{p(x)}(u)
\leq|u|_{p(x)}^{p^+}
\end{equation}
\begin{equation}\label{L5}
|u|_{p(x)}<1\;\;\;\Rightarrow\;\;\;|u|_{p(x)}^{p^+}\leq
\rho_{p(x)}(u)\leq|u|_{p(x)}^{p^-}
\end{equation}
\begin{equation}\label{L6}
|u_n-u|_{p(x)}\rightarrow 0\;\;\;\Leftrightarrow\;\;\;\rho_{p(x)}
(u_n-u)\rightarrow 0.
\end{equation}
Spaces with $p^+ =\infty$ have been studied by Edmunds, Lang and 
Nekvinda
\cite{edm}.

Next, we define $W_0^{1,p(x)}(\Omega)$ as the closure of
$C_0^\infty(\Omega)$ under the norm
$$\|u\|_{p(x)}=|\nabla u|_{p(x)}.$$
The space $(W_0^{1,p(x)}(\Omega),\|\cdot\|_{p(x)})$ is a separable 
and reflexive Banach space. We note that if $q\in C_+(\overline
\Omega)$ and $q(x)<p^\star(x)$ for all $x\in\overline\Omega$ then 
the embedding
$W_0^{1,p(x)}(\Omega)\hookrightarrow L^{q(x)}(\Omega)$
is compact and continuous, where $p^\star(x)=\frac{Np(x)}{N-p(x)}$
if $p(x)<N$ or $p^\star(x)=+\infty$ if $p(x)\geq N$. We refer to
\cite{edm2,edm3,FSZ,FZ1,KR} for further properties of variable exponent
Lebesgue-Sobolev spaces.
\begin{rem}\label{r1}
If $p_1(x)$, $p_2(x)\in C_+(\overline\Omega)$ it is clear that
$m(x)\in C_+(\overline\Omega)$ where $m(x)=\max\{p_1(x),p_2(x)\}$ 
for any $x\in\overline\Omega$. On the other hand since $p_1(x)$,
$p_2(x)\leq m(x)$ for any $x\in\overline\Omega$ it follows that 
$W_0^{1,m(x)}(\Omega)$ is continuously embedded in 
$W_0^{1,p_i(x)}(\Omega)$ for $i\in\{1,2\}$.
\end{rem}

\section{Main results}
In this paper we study Problem \eq{Pr1} if
$f(x,t)=\pm(-\lambda|t|^{m(x)-2}t+|t|^{q(x)-2}t)$, 
where $$m(x):=\max\{p_1(x),p_2(x)\}<q(x)<\left\{
\begin{array}{lll}
&\di \frac{N\cdot m(x)}{N-m(x)}&\di \qquad\mbox{if $m(x)<N$}\\
&\di +\infty&\di \qquad\mbox{if $m(x)\geq N$}\,,\end{array}\right.$$
for any $x\in\overline\Omega$
and all $\lambda>0$.

We first consider the problem
\begin{equation}\label{2}
\left\{\begin{array}{lll}
-{\rm div}((|\nabla u|^{p_1(x)-2}+|\nabla u|^{p_2(x)-2})\nabla u)=
-\lambda |u|^{m(x)-2}u+|u|^{q(x)-2}u, &\mbox{for}& x\in\Omega\\
u=0, &\mbox{for}& x\in\partial\Omega.
\end{array}\right.
\end{equation}
We say that $u\in W_0^{1,m(x)}(\Omega)$ is a 
{\it weak solution} of problem \eq{2} if
$$\int_\Omega(|\nabla u|^{p_1(x)-2}+|\nabla u|^{p_2(x)-2})\nabla u
\nabla v\;dx+\lambda\int_\Omega |u|^{m(x)-2}uv\;dx-\int_\Omega
|u|^{q(x)-2}uv\;dx=0,$$
for all $v\in W_0^{1,m(x)}(\Omega)$.

We prove
\begin{teo}\label{t1}
For every $\lambda>0$ problem \eq{2} has infinitely many weak 
solutions, provided that $2\leq p_i^-$ for $i\in\{1,2\}$, $m^+<q^-$ 
and $q^+<\frac{N\cdot m^-}{N-m^-}$.
\end{teo}

Next, we study the problem
\begin{equation}\label{P3}
\left\{\begin{array}{lll}
-{\rm div}((|\nabla u|^{p_1(x)-2}+|\nabla u|^{p_2(x)-2})\nabla u)=
\lambda |u|^{m(x)-2}u-|u|^{q(x)-2}u, &\mbox{for}& x\in\Omega\\
u=0, &\mbox{for}& x\in\partial\Omega.
\end{array}\right.
\end{equation}

We say that $u\in W_0^{1,m(x)}(\Omega)$ is a 
{\it weak solution} of problem \eq{P3} if
$$\int_\Omega(|\nabla u|^{p_1(x)-2}+|\nabla u|^{p_2(x)-2})\nabla u
\nabla v\;dx-\lambda\int_\Omega |u|^{m(x)-2}uv\;dx+\int_\Omega
|u|^{q(x)-2}uv\;dx=0,$$
for all $v\in W_0^{1,m(x)}(\Omega)$.

We prove
\begin{teo}\label{t2}
There exists $\lambda^\star>0$ such that for any $\lambda\geq
\lambda^\star$ problem \eq{P3} has a nontrivial weak solution,
provided that $m^+<q^-$ and $q^+<\frac{N\cdot m^-}{N-m^-}$.
\end{teo}

\section{Proof of Theorem \ref{t1}}
The key argument in the proof of Theorem \ref{t1} is the following
$\ZZ_2$-symmetric version (for even functionals) of the Mountain 
Pass Lemma (see Theorem 9.12 in \cite{Rab}):

\begin{teo}\label{mpl}  Let $X$ be an infinite 
dimensional real Banach space and let $I\in C^1(X,\RR)$ be even, 
satisfying the Palais-Smale condition (that is, any sequence $\{x_n\}
\subset X$ such that $\{I(x_n)\}$ is bounded and $I^{'}(x_n)
\rightarrow c$ in $X^\star$ has a convergent subsequence) and 
$I(0)=0$. Suppose that
\smallskip

\noindent (I1) There exist two constants $\rho$, $a>0$ such that
$I(x)\geq a$ if $\|x\|=\rho.$
\smallskip

\noindent (I2) For each finite dimensional subspace $X_1\subset X$,
the set $\{x\in X_1;\;I(x)\geq 0\}$ is bounded.
\smallskip

Then $I$ has an unbounded sequence of critical values.
\end{teo} 

Let $E$ denote the generalized Sobolev space $W_0^{1,m(x)}(\Omega)$.

The energy functional corresponding to problem \eq{2} is defined by
$J_\lambda:E\rightarrow\RR$,
$$J_\lambda(u)=\int_\Omega\frac{1}{p_1(x)}|\nabla u|^{p_1(x)}\;dx+
\int_\Omega\frac{1}{p_2(x)}|\nabla u|^{p_2(x)}\;dx+\lambda
\int_\Omega\frac{1}{m(x)}|u|^{m(x)}\;dx-\int_\Omega\frac{1}{q(x)}
|u|^{q(x)}\;dx.$$
A simple calculation based on Remark \ref{r1}, relations \eq{L4}
and \eq{L5} and the compact embedding of $E$ into $L^{s(x)}(\Omega)$
for all $s\in C_+(\overline\Omega)$ with $s(x)<m^\star(x)$ on
$\overline\Omega$ shows that $J_\lambda$ is well-defined on $E$
and $J_\lambda\in C^1(E,\RR)$ with the derivative given by
$$\langle J_\lambda^{'}(u),v\rangle=\int_\Omega(|\nabla u|^{p_1(x)-2}
+|\nabla u|^{p_2(x)-2})\nabla u\nabla v\;dx+\lambda\int_\Omega 
|u|^{m(x)-2}uv\;dx-\int_\Omega|u|^{q(x)-2}uv\;dx,$$
for any $u$, $v\in E$. Thus the weak solutions of \eq{2} are exactly 
the critical points of $J_\lambda$.

\begin{lemma}\label{le1}
There exist $\eta>0$ and $\alpha>0$ such that 
$J_\lambda(u)\geq\alpha>0$ for any $u\in E$ with $\|u\|_{m(x)}=
\eta$.
\end{lemma}
\proof
We first point out that since $m(x)=\max\{p_1(x),p_2(x)\}$ for
any $x\in\overline\Omega$ then
\begin{equation}\label{L8}
|\nabla u(x)|^{p_1(x)}+|\nabla u(x)|^{p_2(x)}\geq|\nabla 
u(x)|^{m(x)},\;\;\;\forall x\in\overline\Omega.
\end{equation}
On the other hand, we have
\begin{equation}\label{L9}
|u(x)|^{q^-}+|u(x)|^{q^+}\geq|u(x)|^{q(x)},\;\;\;\forall x\in\overline
\Omega.
\end{equation}
Using \eq{L8} and \eq{L9} we deduce that
\begin{equation}\label{L10}
\begin{array}{lll}
J_\lambda(u)&\geq&\di\frac{1}{\max\{p_1^+,p_2^+\}}\cdot\di\int_\Omega
|\nabla u|^{m(x)}\;dx-\di\frac{1}{q^-}\cdot\left(\di\int_\Omega|u|^
{q^-}\;dx+\di\int_\Omega|u|^{q^+}\;dx\right)\\
&\geq&\di\frac{1}{m^+}\cdot\di\int_\Omega|\nabla u|^
{m(x)}\;dx-\di\frac{1}{q^-}\cdot\left(\di\int_\Omega|u|^{q^-}\;dx
+\di\int_\Omega|u|^{q^+}\;dx\right)\,,
\end{array}
\end{equation}
for any $u\in E$.

Since $m^+<q^-\leq q^+<m^\star(x)$ for any $x\in\overline\Omega$ and 
$E$ is continuously embedded in $L^{q^-}(\Omega)$ and in 
$L^{q^+}(\Omega)$ it follows that there exist two
positive constants $C_1$ and $C_2$ such that
\begin{equation}\label{L11}
\|u\|_{m(x)}\geq C_1\cdot|u|_{q^+},\;\;\;\|u\|_{m(x)}\geq C_2\cdot
|u|_{q^-},\;\;\;\forall u\in E.
\end{equation}
Assume that $u\in E$ and 
$\|u\|_{m(x)}<1$. Thus, by \eq{L5}, 
\begin{equation}\label{L12}
\int_\Omega|\nabla u|^{m(x)}\;dx\geq\|u\|_{m(x)}^{m^+}.
\end{equation} 
Relations \eq{L10}, \eq{L11} and \eq{L12} yield
\begin{eqnarray*}
J_\lambda(u)&\geq&\frac{1}{m^+}\cdot\|u\|_{m(x)}^{m^+}-\frac{1}{q^-}
\cdot\left[\left(\frac{1}{C_1}\cdot\|u\|_{m(x)}\right)^{q^+}+
\left(\frac{1}{C_2}\cdot\|u\|_{m(x)}\right)^{q^-}\right]\\
&=&(\beta-\gamma\cdot\|u\|_{m(x)}^{q^+-m^+}-\delta\cdot\|u\|_{m(x)}^
{q^--m^+})\cdot\|u\|_{m(x)}^{m^+}
\end{eqnarray*}
for any $u\in E$ with $\|u\|_{m(x)}<1$, where $\beta$, $\gamma$ and
$\delta$ are positive constants.

We remark that the function $g:[0,1]\rightarrow\RR$ defined by
$$g(t)=\beta-\gamma\cdot t^{q^+-m^+}-\delta\cdot t^{q^--m^+}$$
is positive in a neighborhood of the origin. We conclude that Lemma
\ref{le1} holds true.  \endproof

\begin{lemma}\label{le2}
Let $E_1$ be a finite dimensional subspace of $E$. Then the set $S=\{
u\in E_1;\;J_\lambda(u)\geq 0\}$ is bounded.
\end{lemma}

\proof
In order to prove Lemma \ref{le2},  we first show that
\begin{equation}\label{SS1}
\int_\Omega\frac{1}{p_1(x)}|\nabla u|^{p_1(x)}\;dx\leq K_1\cdot
(\|u\|_{m(x)}^{p_1^-}+\|u\|_{m(x)}^{p_1^+}),\;\;\;\forall u\in E
\end{equation}
where $K_1$ is a positive constant.

Indeed, using relations \eq{L4} and \eq{L5} we have
\begin{equation}\label{SS2}
\int_\Omega|\nabla u|^{p_1(x)}\;dx\leq|\nabla u|_{p_1(x)}^{p_1^-}+
|\nabla u|_{p_1(x)}^{p_1^+}=\|u\|_{p_1(x)}^{p_1^-}+
\|u\|_{p_1(x)}^{p_1^+},\;\;\;\forall u\in E.
\end{equation}
On the other hand, Remark \ref{r1} implies that there exists a 
positive constant $K_0$ such that
\begin{equation}\label{SS3}
\|u\|_{p_1(x)}\leq K_0\cdot \|u\|_{m(x)},\;\;\;\forall u\in E.
\end{equation}
Inequalities \eq{SS2} and \eq{SS3} yield
$$\int_\Omega|\nabla u|^{p_1(x)}\;dx\leq
(K_0\cdot\|u\|_{m(x)})^{p_1^-}+
(K_0\cdot\|u\|_{m(x)})^{p_1^+},\;\;\;\forall u\in E$$
and thus \eq{SS1} holds true.

With similar arguments we deduce that there exists a positive constant
$K_2$ such that
\begin{equation}\label{SS4}
\int_\Omega\frac{1}{p_2(x)}|\nabla u|^{p_2(x)}\;dx\leq K_2\cdot
(\|u\|_{m(x)}^{p_2^-}+\|u\|_{m(x)}^{p_2^+}),\;\;\;\forall u\in E.
\end{equation}
Using again \eq{L4} and \eq{L5} we have
$$\int_\Omega|u|^{m(x)}\;dx\leq|u|_{m(x)}^{m^-}+|u|_{m(x)}^{m^+},
\;\;\;\forall u\in E.$$
Since $E$ is continuously embedded in $L^{m(x)}(\Omega)$,
there exists of a positive constant ${\overline{K}}$ such that
$$|u|_{m(x)}\leq{\overline{K}}\cdot\|u\|_{m(x)},\;\;\;\forall u\in E.$$
The last two inequalities show that for each $\lambda>0$ there exists
a positive constant $K_3(\lambda)$ such that
\begin{equation}\label{SS5}
\lambda\cdot\int_\Omega\frac{1}{m(x)}|\nabla u|^{m(x)}\;dx\leq
K_3(\lambda)\cdot(\|u\|_{m(x)}^{m^-}+\|u\|_{m(x)}^{m^+}),\;\;\;
\forall u\in E.
\end{equation} 
By inequalities \eq{SS1}, \eq{SS4} and \eq{SS5} we get
$$J_\lambda(u)\leq K_1\cdot(\|u\|_{m(x)}^{p_1^-}+\|u\|_{m(x)}^{p_1^+})
+K_2\cdot(\|u\|_{m(x)}^{p_2^-}+\|u\|_{m(x)}^{p_2^+})+K_3(\lambda)
\cdot(\|u\|_{m(x)}^{m^-}+\|u\|_{m(x)}^{m^+})-\frac{1}{q^+}
\int_\Omega|u|^{q(x)}\;dx,$$
for all $u\in E$.

Let $u\in E$ be arbitrary but fixed. We define
$$\Omega_<=\{x\in\Omega;\;|u(x)|<1\},\;\;\;\Omega_\geq=\Omega
\setminus\Omega_<.$$
Therefore
\begin{eqnarray*}
J_\lambda(u)&\leq& K_1\cdot(\|u\|_{m(x)}^{p_1^-}+\|u\|_{m(x)}^{p_1^+})
+K_2\cdot(\|u\|_{m(x)}^{p_2^-}+\|u\|_{m(x)}^{p_2^+})+K_3(\lambda)
\cdot(\|u\|_{m(x)}^{m^-}+\|u\|_{m(x)}^{m^+})-\\
&&\frac{1}{q^+}\int_\Omega|u|^{q(x)}\;dx\\
&\leq&K_1\cdot(\|u\|_{m(x)}^{p_1^-}+\|u\|_{m(x)}^{p_1^+})
+K_2\cdot(\|u\|_{m(x)}^{p_2^-}+\|u\|_{m(x)}^{p_2^+})+K_3(\lambda)
\cdot(\|u\|_{m(x)}^{m^-}+\|u\|_{m(x)}^{m^+})-\\
&&\frac{1}{q^+}\int_{\Omega_\geq}|u|^{q(x)}\;dx\\
&\leq&K_1\cdot(\|u\|_{m(x)}^{p_1^-}+\|u\|_{m(x)}^{p_1^+})
+K_2\cdot(\|u\|_{m(x)}^{p_2^-}+\|u\|_{m(x)}^{p_2^+})+K_3(\lambda)
\cdot(\|u\|_{m(x)}^{m^-}+\|u\|_{m(x)}^{m^+})-\\
&&\frac{1}{q^+}\int_{\Omega_\geq}|u|^{q^-}\;dx\\
&\leq&K_1\cdot(\|u\|_{m(x)}^{p_1^-}+\|u\|_{m(x)}^{p_1^+})
+K_2\cdot(\|u\|_{m(x)}^{p_2^-}+\|u\|_{m(x)}^{p_2^+})+K_3(\lambda)
\cdot(\|u\|_{m(x)}^{m^-}+\|u\|_{m(x)}^{m^+})-\\
&&\frac{1}{q^+}\int_{\Omega}|u|^{q^-}\;dx+\frac{1}{q^+}
\int_{\Omega_<}|u|^{q^-}\;dx.
\end{eqnarray*}
But there exists a positive constant $K_4$ such that
$$\frac{1}{q^+}\int_{\Omega_<}|u|^{q^-}\leq K_4,\;\;\;\forall 
u\in E.$$
Hence
\begin{eqnarray*}
J_\lambda(u)&\leq&K_1\cdot(\|u\|_{m(x)}^{p_1^-}+\|u\|_{m(x)}^{p_1^+})
+K_2\cdot(\|u\|_{m(x)}^{p_2^-}+\|u\|_{m(x)}^{p_2^+})+K_3(\lambda)
\cdot(\|u\|_{m(x)}^{m^-}+\|u\|_{m(x)}^{m^+})-\\
&&\frac{1}{q^+}\int_\Omega|u|^{q^-}\;dx+K_4,\;\;\;\forall u\in E.
\end{eqnarray*}
The functional $|\,\cdot\,|_{q^-}:E\rightarrow\RR$ defined by
$$|u|_{q^-}=\left(\int_\Omega|u|^{q^-}\;dx\right)^{1/{q^-}}$$
is a norm in $E$. In the finite dimensional subspace $E_1$ the norms
$|\,\cdot\,|_{q^-}$ and $\|\,\cdot\,\|_{m(x)}$ are equivalent, so there exists a
positive constant $K=K(E_1)$ such that
$$\|u\|_{m(x)}\leq K\cdot|u|_{q^-},\;\;\;\forall u\in E_1.$$
As a consequence we have that there exists a positive constant 
$K_5$ such that
\begin{eqnarray*}
J_\lambda(u)&\leq&K_1\cdot(\|u\|_{m(x)}^{p_1^-}+\|u\|_{m(x)}^{p_1^+})
+K_2\cdot(\|u\|_{m(x)}^{p_2^-}+\|u\|_{m(x)}^{p_2^+})+K_3(\lambda)
\cdot(\|u\|_{m(x)}^{m^-}+\|u\|_{m(x)}^{m^+})-\\
&&K_5\cdot\|u\|_{m(x)}^{q^-}+K_4,\;\;\;\forall u\in E_1.
\end{eqnarray*}
Hence
\begin{eqnarray*}
K_1\cdot(\|u\|_{m(x)}^{p_1^-}+\|u\|_{m(x)}^{p_1^+})
&+&K_2\cdot(\|u\|_{m(x)}^{p_2^-}+\|u\|_{m(x)}^{p_2^+})+K_3(\lambda)
\cdot(\|u\|_{m(x)}^{m^-}+\|u\|_{m(x)}^{m^+})\\
&-&K_5\cdot\|u\|_{m(x)}^{q^-}+K_4\geq 0,\;\;\;\forall u\in S
\end{eqnarray*}
and since $q^->m^+$ we conclude that $S$ is bounded
in $E$.
The proof of Lemma \ref{le2} is complete.
\endproof

\begin{lemma}\label{le3}
Assume that $\{u_n\}\subset E$ is a sequence which satisfies the properties:
\begin{equation}\label{L13}
|J_\lambda(u_n)|<M
\end{equation}
\begin{equation}\label{L14}
J_\lambda^{'}(u_n)\rightarrow 0\;\;\;{\rm as}\;n\rightarrow\infty
\end{equation}
where $M$ is a positive constant. Then $\{u_n\}$ possesses a 
convergent subsequence.
\end{lemma}

\proof
First, we show that $\{u_n\}$ is bounded in $E$. Assume by 
contradiction the contrary. Then, passing eventually at a subsequence,
still denoted by $\{u_n\}$, we may assume that $\|u_n\|_{m(x)}
\rightarrow\infty$ as $n\rightarrow\infty$. Thus we may consider that
$\|u_n\|_{m(x)}>1$ for any integer $n$.

By \eq{L14} we deduce that there exists $N_1>0$ such that for any
$n>N_1$ we have
$$\|J_\lambda^{'}(u_n)\|\leq 1.$$
On the other hand, for any $n>N_1$ fixed, the application
$$E\ni v\rightarrow\langle J_\lambda^{'}(u_n),v\rangle$$
is linear and continuous.
The above information yields 
$$|\langle J_\lambda^{'}(u_n),v\rangle|\leq
\|J_\lambda^{'}(u_n)\|\cdot\|v\|_{m(x)}\leq\|v\|_{m(x)},\;\;\;\forall
v\in E,\;n>N_1.$$
Setting $v=u_n$ we have
$$-\|u_n\|_{m(x)}\leq\int_\Omega|\nabla u_n|^{p_1(x)}\;dx+\int_\Omega
|\nabla u_n|^{p_2(x)}\;dx+\lambda\int_\Omega|u_n|^{m(x)}\;dx-
\int_\Omega|u_n|^{q(x)}\;dx\leq\|u_n\|_{m(x)},$$
for all $n>N_1$. We obtain
\begin{equation}\label{L15}
-\|u_n\|_{m(x)}-\int_\Omega|\nabla u_n|^{p_1(x)}\;dx-\int_\Omega
|\nabla u_n|^{p_2(x)}\;dx-\lambda\int_\Omega|u_n|^{m(x)}\;dx\leq
-\int_\Omega|u_n|^{q(x)}\;dx
\end{equation}
for any $n>N_1$. 

Assuming that $\|u_n\|_{m(x)}>1$, relations \eq{L13}, \eq{L15} and 
\eq{L4} imply
\begin{eqnarray*}
M>J_\lambda(u_n)&\geq&\left(\frac{1}{m^+}-\frac{1}{q^-}\right)
\cdot\int_\Omega(|\nabla u_n|^{p_1(x)}+|\nabla u_n|^{p_2(x)})\;dx\\
&+&\lambda\cdot\left(\frac{1}{m^+}-\frac{1}{q^-}\right)\cdot
\int_\Omega|u_n|^{m(x)}\;dx-\frac{1}{q^-}\cdot\|u_n\|_{m(x)}\\
&\geq&\left(\frac{1}{m^+}-\frac{1}{q^-}\right)\cdot\int_\Omega
|\nabla u_n|^{m(x)}\;dx-\frac{1}{q^-}\|u_n\|_{m(x)}\\
&\geq&\left(\frac{1}{m^+}-\frac{1}{q^-}\right)\cdot\|u_n\|_{m(x)}
^{m^-}-\frac{1}{q^-}\|u_n\|_{m(x)}.
\end{eqnarray*}
Letting $n\rightarrow\infty$ we obtain a contradiction. It follows 
that $\{u_n\}$ is bounded in $E$.

Since $\{u_n\}$ is bounded in $E$, there exist a
subsequence, again denoted by $\{u_n\}$, and $u_0\in E$ such that
$\{u_n\}$ converges weakly to $u_0$ in $E$. Since $E$ is compactly
embedded in $L^{m(x)}(\Omega)$ and in $L^{q(x)}(\Omega)$ it follows
that $\{u_n\}$ converges strongly to $u_0$ in $L^{m(x)}(\Omega)$
and $L^{q(x)}(\Omega)$. The above information and relation \eq{L14}
imply
$$\langle J_\lambda^{'}(u_n)-J_\lambda^{'}(u_0),u_n-u_0\rangle
\rightarrow 0\;\;\;{\rm as}\;n\rightarrow\infty.$$
On the other hand, we have
\begin{equation}\label{L16}
\begin{array}{lll}
\di\int_\Omega(|\nabla u_n|^{p_1(x)-2}\nabla u_n&+&|\nabla u_n|^
{p_2(x)-2}\nabla u_n-|\nabla u_0|^{p_1(x)-2}\nabla u_0-|\nabla u_0|^
{p_2(x)-2}\nabla u_0)\cdot(\nabla u_n-\nabla u_0)\;dx\\
&=&\langle J_\lambda^{'}(u_n)-J_\lambda^{'}(u_0),u_n-u_0\rangle\\
&-&\lambda\cdot\di\int_\Omega(|u_n|^{m(x)-1}u_n-|u_0|^{m(x)-1}u_0)
(u_n-u_0)\;dx\\
&+&\di\int_\Omega(|u_n|^{q(x)-1}u_n-|u_0|^{q(x)-1}u_0)(u_n-u_0)\;dx
\end{array}
\end{equation}
Using the fact that $\{u_n\}$ converges strongly to $u_0$ in $L^{q(x)}
(\Omega)$ and inequality \eq{Hol} we have
$$
\begin{array}{lll}
\left|\di\int_\Omega(|u_n|^{q(x)-1}u_n-|u_0|^{q(x)-1}u_0)
(u_n-u_0)\;dx\right|&\leq&
\left|\di\int_\Omega|u_n|^{q(x)-2}u_n(u_n-u_0)\;dx\right|\\
&+&\left|\di\int_\Omega|u_0|^{q(x)-2}u_0(u_n-u_0)\;dx\right|\\
&\leq&C_3\cdot||u_n|^{q(x)-1}|_{\frac{q(x)}{q(x)-1}}\cdot|u_n-u_0|_
{q(x)}\\
&+&C_4\cdot||u_0|^{q(x)-1}|_{\frac{q(x)}{q(x)-1}}\cdot
|u_n-u_0|_{q(x)}\,,
\end{array}
$$
where $C_3$ and $C_4$ are positive constants. Since $|u_n-u_0|_
{q(x)}\rightarrow 0$ as $n\rightarrow\infty$ we deduce that
\begin{equation}\label{L17}
\lim\limits_{n\rightarrow\infty}
\int_\Omega(|u_n|^{q(x)-1}u_n-|u_0|^{q(x)-1}u_0)
(u_n-u_0)\;dx=0.
\end{equation}
With similar arguments we obtain
\begin{equation}\label{L18}
\lim\limits_{n\rightarrow\infty}
\int_\Omega(|u_n|^{m(x)-1}u_n-|u_0|^{m(x)-1}u_0)
(u_n-u_0)\;dx= 0.
\end{equation} 
By \eq{L16}, \eq{L17} and \eq{L18} we get
\begin{equation}\label{L19}
\begin{array}{lll}
\lim\limits_{n\rightarrow\infty}
\int_\Omega(|\nabla u_n|^{p_1(x)-2}\nabla u_n&+&|\nabla u_n|^
{p_2(x)-2}\nabla u_n-|\nabla u_0|^{p_1(x)-2}\nabla u_0\\
&-&|\nabla u_0|^
{p_2(x)-2}\nabla u_0)\cdot(\nabla u_n-\nabla u_0)\;dx= 0.
\end{array}
\end{equation}
Next, we apply the following elementary inequality (see \cite[Lemma 4.10]{diaz})
\begin{equation}\label{L20}
(|\xi|^{r-2}\xi-|\psi|^{r-2}\psi)\cdot(\xi-\psi)\geq C\,|\xi-\psi|^r,\;\;\;\forall r\geq 
2,\;\xi,\psi\in\RR^N.
\end{equation}
Relations \eq{L19} and \eq{L20} yield
$$\lim\limits_{n\rightarrow\infty}
\int_\Omega|\nabla u_n-\nabla u_0|^{p_1(x)}\;dx+\int_\Omega
|\nabla u_n-\nabla u_0|^{p_2(x)}\;dx=0$$
or using relation \eq{L8} we get
$$\lim\limits_{n\rightarrow\infty}
\int_\Omega|\nabla u_n-\nabla u_0|^{m(x)}\;dx=0.$$
That fact and relation \eq{L6} imply $\|u_n-u_0\|_{m(x)}\rightarrow 
0$ as $n\rightarrow\infty$. 
The proof of Lemma \ref{le3} is complete.  \endproof

{\sc Proof of Theorem \ref{t1} completed.}
It is clear that the functional $J_\lambda$ is even and verifies
$J_\lambda(0)=0$. Lemma \ref{le3} implies that $J_\lambda$ 
satisfies the Palais-Smale condition. On the other hand, 
Lemmas \ref{le1} and \ref{le2} show that conditions (I1) and (I2)
are satisfied. Applying Theorem \ref{mpl} to the
functional $J_\lambda$ we conclude that equation \eq{2} has 
infinitely many weak solutions in $E$.
The proof of Theorem \ref{t1} is complete.\endproof

\section{Proof of Theorem \ref{t2}}
Define the energy functional associated to
Problem \eq{P3} by $I_\lambda:E\rightarrow\RR$,
$$I_\lambda(u)=\int_\Omega\frac{1}{p_1(x)}|\nabla u|^{p_1(x)}\;dx+
\int_\Omega\frac{1}{p_2(x)}|\nabla u|^{p_2(x)}\;dx-\lambda
\int_\Omega\frac{1}{m(x)}|u|^{m(x)}\;dx+\int_\Omega\frac{1}{q(x)}
|u|^{q(x)}\;dx.$$
The same arguments as those used in the case of functional 
$J_\lambda$ show that $I_\lambda$ is well-defined on $E$
and $I_\lambda\in C^1(E,\RR)$ with the derivative given by
$$\langle I_\lambda^{'}(u),v\rangle=\int_\Omega(|\nabla u|^{p_1(x)-2}
+|\nabla u|^{p_2(x)-2})\nabla u\nabla v\;dx-\lambda\int_\Omega 
|u|^{m(x)-2}uv\;dx+\int_\Omega|u|^{q(x)-2}uv\;dx,$$
for any $u$, $v\in E$. We obtain that the weak solutions of \eq{P3} 
are the critical points of $I_\lambda$.

This time our idea is to show that $I_\lambda$ possesses a nontrivial
global minimum point in $E$. With that end in view we start by 
proving two auxiliary results.
\begin{lemma}\label{PL1}
The functional $I_\lambda$ is coercive on $E$.
\end{lemma}
\proof
In order to prove Lemma \ref{PL1} we first show that for any $a$,
$b>0$ and $0<k<l$ the following inequality holds 
\begin{equation}\label{stea1}
a\cdot t^k-b\cdot t^l\leq a\cdot\left(\frac{a}{b}\right)^{k/(l-k)},
\;\;\;\forall\;t\geq 0.
\end{equation}
Indeed, since the function
$$[0,\infty)\ni t\rightarrow t^\theta$$
is increasing for any $\theta>0$ it follows that
$$a-b\cdot t^{l-k}<0,\;\;\;\forall\;t>\left(\frac{a}{b}\right)^
{1/(l-k)},$$
and
$$t^k\cdot(a-b\cdot t^{l-k})\leq a\cdot t^k< a\cdot\left(\frac{a}{b}
\right)^{k/(l-k)},\;\;\;\forall\;t\in\left[0,\left(\frac{a}{b}
\right)^{1/(l-k)}\right].$$
The above two inequalities show that \eq{stea1} holds true.

Using \eq{stea1} we deduce that for any $x\in\Omega$ and $u\in E$
we have
\begin{eqnarray*}
\frac{\lambda}{m^-}|u(x)|^{m(x)}-\frac{1}{q^+}|u(x)|^{q(x)}
&\leq&\frac{\lambda}
{m^-}\left[\frac{\lambda\cdot q^+}{m^-}\right]^{m(x)/(q(x)-m(x))}\\
&\leq&\frac{\lambda}{m^-}\left[\left(\frac{\lambda\cdot q^+}{m^-}
\right)^{m^+/(q^--m^+)}+\left(\frac{\lambda\cdot q^+}{m^-}
\right)^{m^-/(q^+-m^-)}\right]={\cal C},
\end{eqnarray*}
where ${\cal C}$ is a positive constant independent of $u$ and $x$.
Integrating the above inequality over $\Omega$ we obtain
\begin{equation}\label{stea2}
\frac{\lambda}{m^-}\int_\Omega|u|^{m(x)}\;dx-\frac{1}{q^+}
\int_\Omega|u|^{q(x)}\;dx\leq{\cal D}
\end{equation}
where ${\cal D}$ is a positive constant independent of $u$.

Using inequalities \eq{L8} and \eq{stea2} we obtain that for any
$u\in E$ with $\|u\|_{m(x)}>1$ we have
\begin{eqnarray*}
I_\lambda(u)&\geq&\frac{1}{m^+}\int_\Omega|\nabla u|^{m(x)}\;dx-
\frac{\lambda}{m^-}\int_\Omega|u|^{m(x)}\;dx+\frac{1}{q^+}
\int_\Omega|u|^{q(x)}\;dx\\
&\geq&\frac{1}{m^+}\|u\|_{m(x)}^{m^-}-\left(\frac{\lambda}{m^-}
\int_\Omega|u|^{m(x)}\;dx-\frac{1}{q^+}\int_\Omega|u|^{q(x)}\;dx
\right)\\
&\geq&\frac{1}{m^+}\|u\|_{m(x)}^{m^-}-{\cal D}.
\end{eqnarray*}
Thus $I_\lambda$ is coercive and the proof of Lemma \ref{PL1} is
complete.  \endproof
\begin{lemma}\label{PL2}
The functional $I_\lambda$ is weakly lower semicontinuous.
\end{lemma}
\proof
In a first instance we prove that the functionals $\Lambda_i:E
\rightarrow\RR$,
$$\Lambda_i(u)=\int_\Omega\frac{1}{p_i(x)}|\nabla u|^{p_i(x)}\;dx,
\;\;\;\forall\;i\in\{1,2\}$$
are convex. Indeed, since the function 
$$[0,\infty)\ni t\rightarrow t^\theta$$
is convex for any $\theta>1$, we deduce that for each $x\in\Omega$
fixed it holds that
$$\left|\frac{\xi+\psi}{2}\right|^{p_i(x)}\leq
\left|\frac{|\xi|+|\psi|}{2}\right|^{p_i(x)}\leq\frac{1}{2}
|\xi|^{p_i(x)}+\frac{1}{2}|\psi|^{p_i(x)},\;\;\;\forall\xi,\psi
\in\RR^N,\;i\in\{1,2\}.$$
Using the above inequality we deduce that
$$\left|\frac{\nabla u+\nabla v}{2}\right|^{p_i(x)}\leq
\frac{1}{2}|\nabla u|^{p_i(x)}+\frac{1}{2}|\nabla v|^{p_i(x)},
\;\;\;\forall u,v\in E,\;x\in\Omega,\;i\in\{1,2\}.$$
Multiplying with $\frac{1}{p_i(x)}$ and integrating over $\Omega$
we obtain
$$\Lambda_i\left(\frac{u+v}{2}\right)\leq\frac{1}{2}\Lambda_i(u)+
\frac{1}{2}\Lambda_i(v),\;\;\;\forall u,v\in E,\;i\in\{1,2\}.$$
Thus $\Lambda_1$ and $\Lambda_2$ are convex. It follows that 
$\Lambda_1+\Lambda_2$ is convex 

Next, we show that the functional $\Lambda_1+\Lambda_2$ is weakly
lower semicontinuous on $E$. Taking into account that $\Lambda_1+
\Lambda_2$ is convex, by Corollary III.8 in \cite{B} it is enough 
to show that $\Lambda_1+\Lambda_2$ is strongly lower semicontinuous 
on $E$. We fix $u\in E$ and $\epsilon>0$. Let $v\in E$ be arbitrary. 
Since $\Lambda_1+\Lambda_2$ is convex and inequality \eq{Hol} holds 
true we have
\begin{eqnarray*}
\Lambda_1(v)+\Lambda_2(v)&\geq&\Lambda_1(u)+\Lambda_2(u)+\langle
\Lambda_1^{'}(u)+\Lambda_2^{'}(u),v-u\rangle\\
&\geq&\Lambda_1(u)+\Lambda_2(u)-\int_\Omega|\nabla u|^{p_1(x)-1}
|\nabla(v-u)|\;dx--\int_\Omega|\nabla u|^{p_2(x)-1}
|\nabla(v-u)|\;dx\\
&\geq&\Lambda_1(u)+\Lambda_2(u)-D_1\cdot||\nabla u|^{p_1(x)-1}|_
{\frac{p_1(x)}{p_1(x)-1}}\cdot|\nabla(u-v)|_{p_1(x)}-\\
&&D_2\cdot||\nabla u|^{p_2(x)-1}|_{\frac{p_2(x)}{p_2(x)-1}}\cdot
|\nabla(u-v)|_{p_2(x)}\\
&\geq&\Lambda_1(u)+\Lambda_2(u)-D_3\cdot\|u-v\|_{m(x)}\\
&\geq&\Lambda_1(u)+\Lambda_2(u)-\epsilon
\end{eqnarray*}
for all $v\in E$ with $\|u-v\|_{m(x)}<\epsilon/[||\nabla u|^
{p_1(x)-1}|_{\frac{p_1(x)}{p_1(x)-1}}+||\nabla u|^{p_2(x)-1}|_
{\frac{p_2(x)}{p_2(x)-1}}]$, where $D_1$, $D_2$ and $D_3$ are
positive constants. It follows that $\Lambda_1+\Lambda_2$
is strongly lower semicontinuous and since it is convex we obtain 
that $\Lambda_1+\Lambda_2$ is weakly lower semicontinuous.

Finally, we remark that if $\{u_n\}\subset E$ is a sequence which
converges weakly to $u$ in $E$ then $\{u_n\}$ converges strongly 
to $u$ in $L^{m(x)}(\Omega)$ and $L^{q(x)}(\Omega)$. Thus, 
$I_\lambda$ is weakly lower semicontinuous. The proof of Lemma
\ref{PL2} is complete.   \endproof
{\sc Proof of Theorem \ref{t2}.} By Lemmas \ref{PL1} and \ref{PL2} 
we deduce that $I_\lambda$ is coercive and weakly lower 
semicontinuous on $E$. Then Theorem 1.2 in \cite{S} implies that 
there exists $u_\lambda\in E$ a global minimizer of $I_\lambda$ and 
thus a weak solution of problem \eq{P3}. 

We show that $u_\lambda$ is not trivial for $\lambda$ large enough.
Indeed, letting $t_0>1$ be a fixed real and $\Omega_1$ be an open 
subset of $\Omega$ with $|\Omega_1|>0$ we deduce that there exists
$u_0\in C_0^\infty(\Omega)\subset E$ such that $u_0(x)=t_0$ for
any $x\in\overline\Omega_1$ and $0\leq u_0(x)\leq t_0$ in $\Omega
\setminus\Omega_1$. We have
\begin{eqnarray*}
I_\lambda(u_0)&=&\int_\Omega\frac{1}{p_1(x)}|\nabla u_0|^{p_1(x)}\;dx
+\int_\Omega\frac{1}{p_2(x)}|\nabla u_0|^{p_2(x)}\;dx-\lambda
\int_\Omega\frac{1}{m(x)}|u_0|^{m(x)}\;dx+\int_\Omega\frac{1}{q(x)}
|u_0|^{q(x)}\;dx\\
&\leq&L-\frac{\lambda}{m^+}\int_{\Omega_1}|u_0|^{m(x)}\;dx\\
&\leq&L-\frac{\lambda}{m^+}\cdot t_0^{m^-}\cdot|\Omega_1|
\end{eqnarray*}
where $L$ is a positive constant.
Thus, there exists $\lambda^\star>0$ such that $I_\lambda(u_0)<0$
for any $\lambda\in[\lambda^\star,\infty)$. It follows that
$I_\lambda(u_\lambda)<0$ for any $\lambda\geq\lambda^\star$ and
thus $u_\lambda$ is a nontrivial weak solution of problem \eq{P3}
for $\lambda$ large enough. The proof of Theorem \ref{t2} is 
complete. \endproof


\begin{thebibliography}{99}
{\footnotesize
\bibitem{AM} E. Acerbi and G. Mingione, Regularity results for a
class of functionals with nonstandard growth, {\it Arch. Rational
Mech. Anal.} {\bf 156} (2001), 121-140.

\bibitem{B} H. Brezis, {\it Analyse fonctionnelle: th\'eorie et
applications}, Masson, Paris, 1992.

\bibitem{CF} J. Chabrowski and Y. Fu, Existence of solutions for
$p(x)$-Laplacian problems on bounded domains, {\it J. Math. Anal.
Appl.}, in press (doi:10.1016/j.jmaa.2004.10.028).

\bibitem{diaz} J. I. D\'{\i}az, {\it Nonlinear Partial Differential Equations and Free Boundaries. 
Elliptic Equations}, Research Notes
in Mathematics, 106, Pitman, Boston, London, Melbourne, 1986.

\bibitem{D} L. Diening, {\it Theorical and numerical results for
electrorheological fluids}, Ph.D. thesis, University of Frieburg,
Germany, 2002.

\bibitem{dkn} P. Drabek , A. Kufner, and F.  Nicolosi, {\it
Quasilinear Elliptic Equations with Degenerations and Singularities},  
Gruyter Series in Nonlinear Analysis and Applications, Vol. 5, Walter de Gruyter \& Co., Berlin, 
1997.

\bibitem{edm} D. E. Edmunds, J. Lang, and A. Nekvinda, On $L^{p(x)}$
norms, {\it Proc. Roy. Soc. London Ser.~A} {\bf 455} (1999), 219-225.

\bibitem{edm2} D. E. Edmunds and J. R\'akosn\'{\i}k, Density of smooth
functions in $W^{k,p(x)}(\Omega)$, {\it Proc. Roy. Soc. London Ser.~A}
{\bf 437} (1992), 229-236.

\bibitem{edm3} D. E. Edmunds and J. R\'akosn\'{\i}k, Sobolev embedding
with variable exponent, {\it Studia Math.} {\bf 143} (2000), 267-293.

\bibitem{FSZ} X. Fan, J. Shen and D. Zhao, Sobolev Embedding 
Theorems for Spaces $W^{k,p(x)}(\Omega)$, {\it J. Math. Anal. Appl.}
{\bf 262} (2001), 749-760.

\bibitem{FZh} X. L. Fan and Q. H. Zhang, Existence of solutions for
$p(x)$-Laplacian Dirichlet problem, {\it Nonlinear Anal} {\bf 52}
(2003), 1843-1852.

\bibitem{FZZ} X. L. Fan, Q. H. Zhang, and D. Zhao, Eigenvalues of
$p(x)$-Laplacian Dirichlet problem, {\it J. Math. Anal. Appl.}
{\bf 302} (2005), 306-317.

\bibitem{FZ1} X. L. Fan and D. Zhao, On the Spaces $L^{p(x)}(\Omega)$
and $W^{m,p(x)}(\Omega)$, {\it J. Math. Anal. Appl.}, {\bf 263}
(2001), 424-446.

\bibitem{GT} D. Gilbarg and N. S. Trudinger, {\it Elliptic Partial
Differential Equations of Second Order}, Springer, Berlin, 1998.

\bibitem{hal} T. C. Halsey, Electrorheological fluids, {\it Science} 
{\bf 258} (1992), 761-766.

\bibitem{isac} D. Hyers, G. Isac and T. Rassias, {\it Topics in Nonlinear Analysis and 
Applications}, World Scientific Publishing Co.,
Inc., River Edge, NJ, 1997.

\bibitem{KR} O. Kov\'a\v cik and J. R\'akosn\'{\i}k, On spaces 
$L^{p(x)}$ and
$W^{1,p(x)}$, {\it Czechoslovak Math. J.} {\bf 41} (1991), 592-618.

\bibitem{kp} A. Kufner and L.--E. Persson, {\it Weighted Inequalities of Hardy Type}, World 
Scientific Publishing Co., Inc., River Edge, NJ, 2003.

\bibitem{M} J. Musielak, {\it Orlicz Spaces and Modular  Spaces},
Lecture Notes in Mathematics, Vol. 1034, Springer, Berlin, 1983.

\bibitem{nak} H. Nakano, {\it Modulared Semi-ordered Linear Spaces},
Maruzen Co., Ltd., Tokyo, 1950.

\bibitem{orl} W. Orlicz,  \"Uber konjugierte Exponentenfolgen, {\it 
Studia
Math.} {\bf 3} (1931), 200-212.

\bibitem{P} K. Perera, Multiple positive solutions for a class
of quasilinear elliptic boundary-value problems, {\it Electronic
Journal of Differential Equations} {\bf 7} (2003), 1-5.

\bibitem{pfe} C. Pfeiffer, C. Mavroidis, Y. Bar-Cohen, and B. Dolgin,
Electrorheological fluid based force feedback device, in {\it 
Proceedings
of the 1999 SPIE Telemanipulator and Telepresence Technologies VI
Conference (Boston, MA)}, Vol. 3840 (1999), pp. 88-99.

\bibitem{Rab} P. Rabinowitz, Minimax methods in critical point theory
with applications to differential equations, {\it Expository Lectures
from the CBMS Regional Conference held at the University of Miami}, 
American Mathematical Society, Providence, RI, 1984. 

\bibitem{R} M. Ruzicka, {\it Electrorheological Fluids Modeling
and Mathematical Theory}, Springer-Verlag, Berlin, 2002.

\bibitem{sha} I. Sharapudinov, On the topology of the space
$L^{p(t)}([0;1])$,
{\it Matem. Zametki} {\bf 26} (1978),  613-632.

\bibitem{S} M. Struwe, {\it Variational Methods: Applications to
Nonlinear Partial Differential Equations and Hamiltonian Systems},
Springer, Heidelberg, 1996.

\bibitem{tse} I. Tsenov, Generalization of the problem of best 
approximation
of a function in the space $L^s$, {\it Uch. Zap. Dagestan Gos. Univ.} 
{\bf
7} (1961), 25-37.

\bibitem{W} M. Willem, {\it Minimax Theorems}, Birkh\"auser,
Boston, 1996.

\bibitem{WW} W. M. Winslow, Induced Fibration of Suspensions,
{\it J. Appl. Phys.} {\bf 20} (1949), 1137-1140.

\bibitem{Z1} V. Zhikov, Averaging of functionals in the calculus of
variations and elasticity, {\it Math. USSR Izy.} {\bf 29} (1987), 
33-66.

\bibitem{Z2} V. Zhikov, On passing to the limit in nonlinear 
variational problem, {\it Math. Sb.} {\bf 183} (1992), 47-84.


}
\end{thebibliography}
\end{document}